\documentclass[12pt,reqno]{amsart}
\usepackage{amsmath,amssymb,amscd,latexsym,amsthm,mathrsfs}
\usepackage[unicode]{hyperref}
\usepackage{hypbmsec}

\newtheorem{Thm}{Theorem}

\newtheorem{Lem}{Lemma}

\newtheorem{cor}{Corollary}
\theoremstyle{remark}

\renewcommand\d{{\mathrm d}}
\renewcommand\pmod[1]{\;(\operatorname{mod}#1)}

\newcommand\T{\rule{0pt}{2.6ex}}

\newcounter{aa}
\newenvironment{itemspec}[1]{%
\begin{list}{{\rm(\alph{aa})}}{%
\usecounter{aa}\setlength{\itemindent}{24pt}\setlength{\listparindent}{12pt}
      \setlength{\topsep}{10pt}\setlength{\parsep}{2pt}
      \setlength{\leftmargin}{0pt}\setlength{\labelsep}{1ex}}}{\end{list}}

\begin{document}

\title[The Erd\H{o}s--Moser equation revisited]%
{The Erd\H{o}s--Moser equation $1^k+2^k+\dots+(m-1)^k=m^k$ \\
revisited using continued fractions}

\author{Yves Gallot}
\address{12 bis rue Perrey, 31400 Toulouse, France}
\email{galloty@orange.fr}

\author{Pieter Moree}
\address{Max-Planck-Institut f\"ur Mathematik, Vivatsgasse 7, D-53111 Bonn, Germany}
\email{moree@mpim-bonn.mpg.de}

\author{Wadim Zudilin}
\address{School of Mathematical and Physical Sciences,
University of Newcastle, Calla\-ghan NSW 2308, AUSTRALIA}
%\email{wzudilin@gmail.com}
\email{wadim.zudilin@newcastle.edu.au}

\date{\today}

\begin{abstract}
If the equation of the title has an integer solution with $k\ge 2$,
then $m>10^{9.3\cdot 10^6}$. This was the current best result and proved using
a method due to L.~Moser (1953). This approach cannot be improved to
reach the benchmark $m>10^{10^7}$. Here we achieve $m>10^{10^9}$ by
showing that $2k/(2m-3)$ is a convergent of $\log 2$ and making an
extensive continued fraction digits calculation of $(\log 2)/N$, with
$N$~an appropriate integer. This method is very different from that of Moser. Indeed,
our result seems to give one of very few instances where a large scale computation
of a numerical constant has an application.
\end{abstract}


\maketitle
%==================================================

\section{Introduction}
\label{s1}

In this note we are interested in non-trivial integer solutions, that is,
solutions with $k\ge 2$, of the equation
\begin{equation}
\label{EME}
1^k+2^k+\dots+(m-2)^k+(m-1)^k=m^k.
\end{equation}
Conjecturally such solutions do not exist. For $k=1$ one has
clearly the solution $1+2=3$ (and no further ones). From now on
we will assume that $k\ge 2$.
Moser~\cite{Moser} showed in 1953 that if $(m,k)$ is a solution of \eqref{EME}, then
$m>10^{10^6}$ and $k$ is even. His result has since then been improved on.
Butske et al.~\cite{graphs} have shown by computing, rather
than estimating, certain quantities in Moser's original proof that $m>1.485\cdot 10^{\,9\,321\,155}$.
By proceeding along these lines this bound cannot be improved on substantially.
Butske et al.~\cite[p. 411]{graphs} expressed
the hope that new insights will eventually make it possible to reach the more natural benchmark
$10^{10^7}$.

Using that $\Sigma_k(m)=1^k+2^k+\dots+m^k\le\int_1^m t^k \d t$
and $\Sigma_k(m+1)>\int_0^m t^k\d t$ we obtain that $k+1<m<2(k+1)$.
This shows that the ratio $k/m$ is bounded.
By a more elaborate reasoning along these lines Krzysztofek~\cite{K}
obtained that $k+2<m<\frac32(k+1)$. This
implies that $k\ge 4$ and hence
\begin{equation}
\label{ineq}
k+2<m<2k.
\end{equation}
Dividing both sides of \eqref{EME} by $m^k$ one sees that for every integer $m\ge 2$, \eqref{EME}
has precisely one \emph{real} solution~$k$. It is known that
$\lim_{m\to \infty}k/m=\log 2$ and we show here that in fact the behaviour of $k$ as a function
of $m$ can be determined in a much more explicit way (Theorem~\ref{th:As} and Section~\ref{sec:As}).

Moree et al.~\cite{MRU}, using properties of the
Bernoulli numbers and polynomials (an approach initiated in
Urbanowicz \cite{Urbi}), showed that
$N_1=\operatorname{lcm}(1,2,\dots,200)\mid k$. Kellner \cite{Kellner} in 2002
showed that also all primes $200<p<1000$ have to divide~$k$. Actually, Moree et al.~\cite[p.~814]{MRU} proved
a slightly stronger result and on combining this with Kellner's, one obtains that $N_2\mid k$ with
\begin{align*}
N_2
=2^8\cdot 3^5\cdot 5^4\cdot 7^3\cdot 11^2\cdot 13^2\cdot 17^2\cdot 19^2\cdot\prod_{23\le p\le 997}p
>5.7462\cdot 10^{427}.
\end{align*}
For some further references and info on the Erd\H{o}s--Moser equation we refer to the book
by Guy~\cite[D7]{Guy}.

In this note we attack \eqref{EME} using the theory of
continued fractions. This approach was first explored in~1976 by Best and
te Riele \cite{Best} in their attempt to solve the related conjecture of
Erd\H{o}s~\cite{E} that there are infinitely many pairs $(m,k)$ such that
$\Sigma_k(m)\ge m^k$ and $2(m-1)^k<m^k$. In this
context they also gave the following variant of one of their results (without
proof), namely, \eqref{eq:2} with $O(m^{-2})$ replaced with $o(m^{-1})$.
The proof we give here uses the
same circle of ideas as used by Best and te Riele. It seems that after their work continued
fractions in the Erd\H{o}s--Moser context have been completely ignored.
We hope the present paper makes clear that this is unjustified.

\begin{Thm}
\label{th:As}
For integer $m>0$ and \emph{real} $k>0$ satisfying equation~\eqref{EME},
we have the asymptotic expansion
\begin{equation}
\label{eq:2}
k=\log2\biggl(m-\frac32-\frac{c_1}m+O\biggl(\frac 1{m^2}\biggr)\biggr)
\quad\text{as}\; m\to\infty,
\end{equation}
with $c_1=\frac{25}{12}-3\log 2\approx 0.00389\dots$\,.
Moreover, if $m>10^9$ then
\begin{equation}
\frac km=\log2\biggl(1-\frac3{2m}-\frac {C_m}{m^2}\biggr),
\qquad\text{where}\quad 0<C_m<0.004.
\label{E01}
\end{equation}
\end{Thm}

\begin{cor}
\label{cor:1}
If $(m,k)$ is a solution of {\rm\eqref{EME}} with $k\ge 2$, then
$2k/(2m-3)$ is a convergent $p_j/q_j$ of $\log 2$ with $j$ even.
\end{cor}
\begin{cor}
\label{cor:2}
The number of solutions $m\le x$ of \eqref{EME}, as $x$ tends to infinity, is at most $O(\log x)$.
\end{cor}
The equation \eqref{EME} seems to be a sole example of an exponential Diophantine
equation in just two unknowns for which even the finiteness of solutions is not yet established. The best
result in this direction is given by Corollary \ref{cor:2}, which is an immediate consequence of the exponential
growth of $p_j$ as a function of $j$ and Corollary \ref{cor:1}.

Corollary \ref{cor:1} is not the only result which relates convergents to solutions of
Diophantine equations. For example, if $(x_0,y_0)$ is a positive solution to Pell's
equation $x^2-dy^2=\pm 1$, with $d$ a positive square-free integer, then $x_0/y_0$ is
a convergent of the continued fraction expansion of $\sqrt{d}$.
On the other hand, in our situation the number in question, $\log2$, is transcendental and
its continued fraction expansion is expected to be sufficiently `generic'
(unlike that of quadratic irrationals).

Corollary \ref{cor:1} naturally leads us to investigate common factors of $k$ and $2m-3$. This can
be done using the method of Moser, but is not in the literature, as before there was no special
reason for considering $2m-3$.

A key role in this arithmetic study is played by the congruence
\begin{equation}
\label{staudt}
\frac{\sum_{j=1}^{l-1}j^r}y=\begin{cases}
0\pmod{\frac12} & \text{if $r>1$ is odd}; \\
-\sum_{p\mid l, \, p-1\mid r}\frac1p\pmod1 &\text{otherwise}.
\end{cases}
\end{equation}
This identity can be proved using the Von Staudt--Clausen theorem; for
alternative proofs see, e.g., Carlitz \cite{Carlitz} or Moree \cite{Canada}.
Its relevance for the study of \eqref{EME} was first pointed out by Moree~\cite{Oz}.

Given $N\ge 1$, put
$$
{\mathcal{P}}(N)=\{p:p-1\mid N\}\cup\{p:\text{$3$ is a primitive root modulo $p$}\}.
$$
By a classical result of Hooley \cite{Hooley} it follows, assuming the
Generalized Riemann Hypothesis (GRH), that ${\mathcal{P}}(N)$ has
a natural density $A$, with $A=0.3739558136\dots$ the Artin constant, in the set of primes.
If $2k/(2m-3)=p_j/q_j$ is a convergent of $\log 2$ arising in Corollary~\ref{cor:1}, then it can be shown that
$(q_j,6)=1$ and,
if $p\in {\mathcal{P}}(N_2)$ and $p$ divides $q_j$, then
$$
\nu_p(q_j)=\nu_p(3^{p-1}-1)+1\ge 2,
$$
where we write $\nu_p(n)=a$ if $p^a\mid n$ and $p^{a+1}\nmid n$.
All primes $p\le 2017$ are in ${\mathcal{P}}(N_2)$. For $p\ne 3$ we have $\nu_p(3^{p-1}-1)=1$ unless
$3^{p-1}\equiv 1\pmod{p^2}$, that is, $p$~is a Mirimanoff prime. (It is known that the only
Mirimanoff primes $p<10^{14}$ are $11$ and $1006003$.)

The main idea of this paper is, in essence, to make use of the fact that the convergents
$p_j/q_j$ of $\log 2$ have no reason to also satisfy $N_2\mid p_j$. The first piece of information
comes from asymptotic analysis and the latter piece from arithmetic. Analysis and
arithmetic give rise to conditions on the solutions that `do not feel
each other' and this is exploited in our main result:

\begin{Thm}
\label{main}
Let $N\ge 1$ be an arbitrary integer. Let
$$\frac{\log 2}{2N} = [a_0,a_1,a_2,\dots] = a_0 + \cfrac{1}{a_1 + \cfrac{1}{a_2 + \cdots}}$$
be the (regular) continued fraction of
$({\log 2})/(2N)$, with $p_i/q_i = [a_0,a_1,\dots,a_i]$ its $i$-th partial convergent.

Suppose that the integer pair $(m,k)$ with $k\ge 2$ satisfies \eqref{EME}
with $N\mid k$.
Let $j=j(N)$ be the smallest integer such that:
\begin{itemspec}

\itemindent18pt
\item
\label{ca}
$j$ is even;

\item
\label{cb}
$a_{j+1}\ge 180N-2$;

\item
\label{cc}
$(q_j,6)=1$; and

\item
\label{cd}
$\nu_p(q_j)=\nu_p(3^{p-1}-1)+\nu_p(N)+1$
for all primes $p\in {\mathcal{P}}(N)$ dividing $q_j$.
\end{itemspec}
Then $m>q_j/2$.
\end{Thm}

Computing many partial quotients (that is, continued fraction digits)
of $\log 2$ is closely related to computing $\log 2$ with
many digits of accuracy. Indeed, it is a well-known result of Lochs that for a generic number knowing
it accurately up to $n$ decimal digits implies that we can compute about $0.97n$
(where $0.97\approx6(\log 2)(\log 10)/\pi^2$)
continued fraction digits accurately. For example, knowing 1000 decimal digits of~$\pi$ allows one
to compute 968 continued fraction digits.

It seems a hopeless problem to prove anything about ${\mathsf E}(\log q_{j(N)})$, the
expected value of $\log q_{j(N)}$ produced by the result. However, metric theory
of continued fractions offers some hope of proving a non-trivial lower bound for
${\mathsf E}(\log q_{j(N)}(\xi))$,
where we require conditions \ref{ca}, \ref{cb}, \ref{cc} and~\ref{cd} to be satisfied but replace
$(\log 2)/(2N)$ by a `generic' $\xi\in [0,1]\setminus \mathbb Q$.
In this context recall the result of L\'evy~\cite{L} that, for such a $\xi$,
\begin{equation}
\lim_{j\to\infty}\frac{\log q_j(\xi)}j=\frac{\pi^2}{12\log 2}\approx1.18.
\label{Levy}
\end{equation}
The Gauss--Kuz'min statistics asserts that,
for a generic~$\xi$, the probability that a given term
in its continued fraction expansion is at least~$b$,
equals $\log_2(1+1/b)$. This allows one to deal with the case where we only have
condition~\ref{cb}. Likewise a result of Moeckel~\cite{Moeckel},
reproved in a very different way a few years later by Jager and
Liardet \cite{JL}, allows one to deal with the case
where we only focus on condition~\ref{cc}. Their result says
that for a generic $\xi\in[0,1]\setminus \mathbb Q$ we have
$$
\lim_{n\to\infty}\frac{\{1\le m\le n:q_m(\xi)\equiv a\pmod{d}\}}{n}
=\frac d{J(d)}\,\frac{\varphi((a,d))}{(a,d)},
$$
where $\varphi$ denotes Euler's totient function, $J(m)=m^2\prod_{p\mid m}(1-1/p^2)$ Jordan's
totient and $(a,m)$ the greatest common divisor of $a$ and~$m$.
%%% PM: ** 2009-06-30
This result shows that $(q_j,6)=1$ with probability $1/2$
(note that a natural number is coprime to~6 with probability $1/3$).
P.~Liardet communicated
to us that methods of his paper \cite{Liardet} can be used to take into account
both conditions \ref{ca} and~\ref{cc};
also the authors of \cite{Harman} claim that this can be done.
We expect that there is a positive constant $c_1$ such that for a generic
$\xi$ satisfying conditions
\ref{ca}, \ref{cb} and \ref{cc}, we have ${\mathsf E}(\log q_{j(N)}(\xi))\sim c_1N$ as $N$ tends to
infinity. Furthermore, we expect that for a generic $\xi$ satisfying conditions
\ref{ca}, \ref{cb}, \ref{cc} and \ref{cd},
${\mathsf E}(\log q_{j(N)}(\xi))\sim c_2 N\log^{\beta} N$
for some positive constants $c_2$ and $\beta$; condition~\ref{cb}
is responsible for $N$, condition~\ref{cd} for $\log^{\beta}N$, while conditions \ref{ca} and~\ref{cc}
affect $c_2$. We are definitely not experts in metric aspects of number theory,
thus leave this problem to the interested reader acquainted with the subject.
Indeed, we even expect that going beyond computing the expected value of $\log q_{j(N)}(\xi)$ is
possible, and a probability distribution for $\log q_{j(N)}(\xi)$ can be obtained.

Using the above results from metric theory of continued fractions and some heuristics we are led
to believe that roughly speaking we can get
$$
m>10^{257N}
$$
from Theorem~\ref{main}. Being able to compute the convergents of $(\log 2)/(2N)$ arbitrarily far,
we would expect (taking $N=N_2$) to show that $m>10^{10^{400}}$.
With the current computer technology computing sufficiently many convergents is the bottleneck.
Taking this into consideration we would expect to get
$$
m>10^{0.515 r},
$$
from Theorem~\ref{main}, where $r$ is the number of convergents we can compute accurately
and $0.515$ is the base~10 logarithm of L\'evy's constant~\eqref{Levy}.
Note that the fact that $N_2$ has many divisors gives us some flexibility and increases
the likelihood of the heuristics to be applicable. Indeed, our numerical experimenting agrees well
with our heuristic considerations (see Section~\ref{s4}). Early 2009, A.~Yee and R.~Chan \cite{AJY}
reached $r>31\cdot10^9$ for~$\log 2$. On the other hand,
Y.~Kanada and his team~\cite{Kanada} computed $\pi$ to over 1.24 trillion decimal
digits already in 2002, using formulae of the same complexity as those used
for the computation of $\log2$ (see \cite[Chapter~3]{BB} for details).
Thus, given the present computer (im)possibilities,
one could hope to show (with a lot of effort!) that $m>10^{10^{12}}$.

Applying Theorem~\ref{main} with $N=2^8\cdot 3^5\cdot 5^3$
or $N=2^8\cdot 3^5\cdot 5^4$, and invoking the result of
Moree et al.~\cite{MRU} that $N\mid k$, we obtain the following

\begin{Thm}
\label{YVES}
If an integer pair $(m,k)$ with $k\ge 2$ satisfies \eqref{EME}, then
$$m > 2.7139 \cdot 10^{\,1\,667\,658\,416}.$$
\end{Thm}

As an application we can show that $\omega(m-1)\ge 33$,
this improves on the result of Brenton and Vasiliu \cite{BV},
who have shown that $\omega(m-1)\ge 26$, where $\omega$
denotes the number of distinct prime divisors; see Section~\ref{s5.1} for further details.

The fact $N_2\mid k$ naively implies that $k$~is of size $10^{427}$ (at least),
which is much smaller than Moser's $10^{10^6}$. However, in this
paper we show that the fact actually yields that $k>10^{10^9}$ (and likely even
$k>10^{10^{400}}$)\,---\,a modestly {\small small} number dividing~$k$
leads to a huge lower bound for~$k$.
Thus, on revisiting \cite{MRU} after 16~years, its main result is seen
to be far more powerful than the second author thought at that time.

In the three following sections we prove Theorems~\ref{th:As},
\ref{main} and \ref{YVES}, respectively.
Our final Section~\ref{s5} is devoted to discussing some problems
related to the Erd\H os--Moser equation.

\section{Asymptotic dependence of $k$ in terms of $m$}
\label{sec:As}

Our proof of Theorem~\ref{th:As} makes use of the following lemma.
\begin{Lem}
\label{lem:As}
For any real $k>0$, we have
\begin{equation}
(1-y)^k=e^{-ky}\biggl(1-\frac k2y^2-\frac k3y^3
+\frac{k(k-2)}8y^4+\frac{k(5k-6)}{30}y^5+O(y^6)\biggr)
\quad \text{as}\; y\to0.
\label{E02}
\end{equation}
Moreover, for $k>8$ and $0<y<1$, the inequality
\begin{align}
&
e^{-ky}\biggl(1-\frac k2y^2-\frac k3y^3+\frac{k(k-2)}8y^4+\frac{k(5k-6)}{30}y^5-\frac{k^3}6y^6\biggr)
\nonumber\\ &\quad
<(1-y)^k
<e^{-ky}\biggl(1-\frac k2y^2-\frac k3y^3+\frac{k(k-2)}8y^4+\frac{k^2}2y^5\biggr)
\label{E03}
\end{align}
holds.
\end{Lem}

\begin{proof}
As for the asymptotic relation in~\eqref{E02}, we simply develop
the Taylor expansion of $(1-y)^ke^{ky}$ up to~$y^5$.
Unfortunately, estimates coming from the classical forms
for the remainder are not sufficient to derive a sharp
dependence on $k$ as in~\eqref{E03} for the last term. Therefore,
we need more drastic methods to quantify the asymptotics in~\eqref{E02}
when $0<y<1$.

First note that
\begin{align*}
(1-y)e^y
&=1-\sum_{n=2}^\infty\frac{n-1}{n!}y^n
\\
&=1-\frac{y^2}2-\frac{y^3}3-\frac{y^4}8-\frac{y^5}{30}-\dotsb,
\qquad 0<y<1.
\end{align*}
Since all coefficients, starting from $n=2$, in this power series are
negative and their sum is exactly $-1$, for these values of~$y$ we have the inequality
\begin{equation}
1-\frac{y^2}2-\frac{y^3}3-\frac{y^4}8-\frac{y^5}{30}-\frac{y^6}{120}
<(1-y)e^y
<1-\frac{y^2}2-\frac{y^3}3-\frac{y^4}8.
\label{E04}
\end{equation}
The quantities
\begin{equation}
x_1=\frac{y^2}2+\frac{y^3}3+\frac{y^4}8
\quad\text{and}\quad
x_2=\frac{y^2}2+\frac{y^3}3+\frac{y^4}8+\frac{y^5}{30}+\frac{y^6}{120},
\label{E05}
\end{equation}
which appear in~\eqref{E04}, lie between $0$ and~$1$ for $0<y<1$.

Our next ingredient is Gerber's generalization of the Bernoulli inequality~\cite{Gerber} (see
also Alzer \cite{Alzer}). It states that the remainder after $k$ terms of the (possibly divergent) binomial
series for $(1+x)^a$ ($a,x$ real with $-1<x$) has the same sign as the first neglected term.
In particular we have for real $k>2$ and $0<x<1$,
\begin{equation}
(1-x)^k<1-kx+\frac{k(k-1)}2x^2,
\label{E06}
\end{equation}
and for real $k>3$ and $0<x<1$,
\begin{equation}
(1-x)^k>1-kx+\frac{k(k-1)}2x^2-\frac{k(k-1)(k-2)}6x^3.
\label{E07}
\end{equation}

Using the right inequality in~\eqref{E04} and taking $x=x_1$ in~\eqref{E06} we obtain, for $k>2$,
\begin{align}
(1-y)^ke^{ky}
&<1-k\biggl(\frac{y^2}2+\frac{y^3}3+\frac{y^4}8\biggr)
+\frac{k(k-1)}2\biggl(\frac{y^2}2+\frac{y^3}3+\frac{y^4}8\biggr)^2
\nonumber\\
&=1-\frac k2y^2-\frac k3y^3+\frac{k(k-2)}8y^4
\nonumber\\ &\qquad
+k(k-1)y^5\biggl(\frac16+\frac{17}{144}y+\frac1{24}y^2+\frac1{128}y^3\biggr)
\nonumber\\
&<1-\frac k2y^2-\frac k3y^3+\frac{k(k-2)}8y^4
+\frac{385}{1152}k(k-1)y^5
\label{E08}
\end{align}
implying the upper estimate in~\eqref{E03}.
In the same vein, the application of the left identity in~\eqref{E04} and of~\eqref{E07}
with $x=x_2$ results, for $k>3$, in
\begin{align}
(1-y)^ke^{ky}
&>1-\frac k2y^2-\frac k3y^3+\frac{k(k-2)}8y^4+\frac{k(5k-6)}{30}y^5
\nonumber\\ &\qquad
-ky^6\sum_{n=0}^{12}(a_nk^2+b_nk+c_n)y^n,
\label{E09}
\end{align}
where the polynomials $p_n(k)=a_nk^2+b_nk+c_n$, $n=0,1,\dots,12$, all have
positive leading coefficients $a_n$; moreover, $p_n(k)>0$ for $k>3$ and
$n=2,3,\dots,12$, \ $p_1(k)=\frac1{24}k^2-\frac{11}{60}k+\frac{17}{120}>0$
for $k>4$, and $p_0(k)=\frac1{48}k^2-\frac{13}{72}k+\frac{121}{720}>0$
for $k>8$. Using this positivity of the polynomials we can continue
the inequality in~\eqref{E09} for $k>8$ as follows:
\begin{align}
(1-y)^ke^{ky}
&>1-\frac k2y^2-\frac k3y^3+\frac{k(k-2)}8y^4+\frac{k(5k-6)}{30}y^5
\nonumber\\ &\qquad
-ky^6\sum_{n=0}^{12}(a_nk^2+b_nk+c_n)
\nonumber\\
&=1-\frac k2y^2-\frac k3y^3+\frac{k(k-2)}8y^4+\frac{k(5k-6)}{30}y^5
\nonumber\\ &\qquad
-ky^6\biggl(\frac16k^2-\frac{17}{24}k+\frac{11}{20}\biggr),
\label{E10}
\end{align}
from which we deduce the left inequality in~\eqref{E03}, and the lemma
follows.
\end{proof}

\begin{proof}[Proof of Theorem~{\rm\ref{th:As}}]
The original equation~\eqref{EME} is equivalent to
\begin{equation}
\label{eq:1}
1=\sum_{j=1}^{m-1}\biggl(1-\frac jm\biggr)^k.
\end{equation}
Applying to each term on the right-hand side the inequality
from~\eqref{E03} we obtain
\begin{align}
&
S_0-\frac k{2m^2}S_2-\frac k{3m^3}S_3+\frac{k(k-2)}{8m^4}S_4+\frac{k(5k-6)}{30m^5}S_5-\frac{k^3}{6m^6}S_6
\nonumber\\ &\qquad
<\sum_{j=1}^{m-1}\biggl(1-\frac jm\biggr)^k
<S_0-\frac k{2m^2}S_2-\frac k{3m^3}S_3+\frac{k(k-2)}{8m^4}S_4+\frac{k^2}{2m^5}S_5,
\label{E11}
\end{align}
with the notation
$$
S_n=\sum_{j=1}^{m-1}j^ne^{-kj/m}
=\sum_{j=1}^{m-1}j^nz^j\bigg|_{z=e^{-k/m}}.
$$

By (\ref{ineq}) we have $e^{-1}<z<e^{-1/2}$, where $z=e^{-k/m}$, and
hence $1/(1-z)<1/(1-e^{-1/2})<3$, and in the closed-form expression of the sum
$$
S_0=\sum_{j=1}^{m-1}z^j=\frac z{1-z}-\frac{z^m}{1-z},
$$
the second term as well as its $z\frac{\d}{\d z}$-derivatives are bounded:
\begin{gather*}
0<\frac{z^m}{1-z}\bigg|_{z=e^{-k/m}}<3e^{-k}
\quad\text{and}
\\
0<\biggl(\biggl(z\frac{\d}{\d z}\biggr)^n\frac{z^m}{1-z}\biggr)\bigg|_{z=e^{-k/m}}<3^{n+1}m^ne^{-k},
\quad\text{for $n=1,2,\dots$}\,.
\end{gather*}
Therefore, we can write the inequality in~\eqref{E11} as
\begin{align*}
&
S_0'-\frac k{2m^2}S_2'-\frac k{3m^3}S_3'+\frac{k(k-2)}{8m^4}S_4'+\frac{k(5k-6)}{30m^5}S_5'-\frac{k^3}{6m^6}S_6'
-\biggl(\frac{3^3k}2+\frac{3^4k}3+\frac{3^7k^3}6\biggr)e^{-k}
\\ &\qquad
<\sum_{j=1}^{m-1}\biggl(1-\frac jm\biggr)^k
\\ &\qquad
<S_0'-\frac k{2m^2}S_2'-\frac k{3m^3}S_3'+\frac{k(k-2)}{8m^4}S_4'+\frac{k^2}{2m^5}S_5'
+\biggl(3+\frac{3^5k(k-2)}8+\frac{3^6k^2}2\biggr)e^{-k}
\end{align*}
implying
\begin{align}
&
S_0'-\frac k{2m^2}S_2'-\frac k{3m^3}S_3'+\frac{k(k-2)}{8m^4}S_4'+\frac{k(5k-6)}{30m^5}S_5'-\frac{k^3}{6m^6}S_6'
-500k^3e^{-k}
\nonumber\\ &\qquad
<\sum_{j=1}^{m-1}\biggl(1-\frac jm\biggr)^k
<S_0'-\frac k{2m^2}S_2'-\frac k{3m^3}S_3'+\frac{k(k-2)}{8m^4}S_4'+\frac{k^2}{2m^5}S_5'
+500k^2e^{-k},
\label{E12}
\end{align}
where
\begin{align*}
S_n'
&=\sum_{j=1}^\infty j^nz^j\bigg|_{z=e^{-k/m}}
=\biggl(\biggl(z\frac{\d}{\d z}\biggr)^n\frac z{1-z}\biggr)\bigg|_{z=e^{-k/m}}
\\
&=(-1)^n\biggl(\biggl(z\frac{\d}{\d z}\biggr)^n\frac1{z-1}\biggr)\bigg|_{z=e^{k/m}}
\quad\text{for $n=0,1,\dots$};
\end{align*}
in particular,
\begin{gather*}
S_0'=\frac1{z-1},
\quad
S_2'=\frac{z+z^2}{(z-1)^3},
\quad
S_3'=\frac{z+4z^2+z^3}{(z-1)^4},
\quad
S_4'=\frac{z+11z^2+11z^3+z^4}{(z-1)^5},
\\
S_5'=\frac{z+26z^2+66z^3+26z^4+z^5}{(z-1)^6},
\quad
S_6'=\frac{z+57z^2+302z^3+302z^4+57z^5+z^6}{(z-1)^7}
\end{gather*}
with $z=e^{k/m}$. Since
$500k^3e^{-k}<(2k)^{-3}<m^{-3}$
for $k>m/2>30$, using our equation \eqref{eq:1} we can write
the estimates~\eqref{E12} as
\begin{align}
&
\frac{k(5k-6)}{30m^5}S_5'-\frac{k^3}{6m^6}S_6'-\frac1{m^3}
\nonumber\\ &\qquad
<1-S_0'+\frac k{2m^2}S_2'+\frac k{3m^3}S_3'-\frac{k(k-2)}{8m^4}S_4'
<\frac{k^2}{2m^5}S_5'+\frac1{m^3}.
\label{E13}
\end{align}
Noting that $e^{1/2}<z=e^{k/m}<e$, we find
\begin{align*}
0<S_5'&<\frac{e+26e^2+66e^3+26e^4+e^5}{(e^{1/2}-1)^6}
<41438,
\\
0<S_6'&<\frac{e+57e^2+302e^3+302e^4+57e^5+e^6}{(e^{1/2}-1)^7}
<658544,
\end{align*}
we continue~\eqref{E13} as follows:
\begin{multline}
\biggl|1-\frac1{z-1}+\frac k{2m^2}\frac{z+z^2}{(z-1)^3}
+\frac k{3m^3}\frac{z+4z^2+z^3}{(z-1)^4}
\\
-\frac{k(k-2)}{8m^4}\frac{z+11z^2+11z^3+z^4}{(z-1)^5}\biggr|
<\frac{110000}{m^3},
\label{E14}
\end{multline}
where $z=e^{k/m}$.

We already know that $k/m$ is bounded as $m\to\infty$;
making the ansatz $k/m=c+O(1/m)$, hence $z=e^{k/m}=e^c+O(1/m)$,
we find from~\eqref{E14} that
$$
1-\frac1{e^c-1}=O\biggl(\frac 1m\biggr) \quad\text{as}\; m\to\infty,
$$
hence $e^c=2$ and $c=\log2$. Now we take
$$
\frac km=\log2+\frac am+\frac b{m^2}+O\biggl(\frac 1{m^3}\biggr) \quad\text{as}\; m\to\infty,
$$
hence
$$
z=e^{k/m}=2+\frac{2a}m+\frac{a^2+2b}{m^2}+O\biggl(\frac 1{m^3}\biggr) \quad\text{as}\; m\to\infty.
$$
Substituting these formulas into~\eqref{E14} results in
\begin{align*}
O\biggl(\frac 1{m^{3}}\biggr)
&=1-\frac1{1+2a/m+(a^2+2b)/m^2+O(m^{-3})}
\\ &\qquad
+\frac{\log2+a/m+O(m^{-2})}{2m}\,\frac{6+10a/m+O(m^{-2})}{1+6a/m+O(m^{-2})}
\\ &\qquad
+\frac{\log2+O(m^{-1})}{3m^2}\,\frac{26+O(m^{-1})}{1+O(m^{-1})}
\\ &\qquad
-\frac{\log^22+O(m^{-1})}{8m^2}\,\frac{150+O(m^{-1})}{1+O(m^{-1})}
+O\biggl(\frac 1{m^3}\biggr)
\\
&=\frac{2a+3\log2}m
-\frac{3a^2-3a+13a\log2-2b+\frac{75}4\log^22-\frac{26}3\log2}{m^2}
+O\biggl(\frac 1{m^3}\biggr),
\end{align*}
hence $a=-\frac32\log2$, $b=(3\log2-\frac{25}{12})\log2$ and, finally,
we get the asymptotic formula~\eqref{eq:2}.

To quantify this asymptotic expansion, we introduce the function
\begin{align*}
f_m(C)
&=\biggl(1-\frac1{z-1}
+\frac\lambda{2m}\,\frac{z+z^2}{(z-1)^3}
+\frac\lambda{3m^2}\,\frac{z+4z^2+z^3}{(z-1)^4}
\\ &\qquad
-\frac{\lambda(\lambda-2/m)}{8m^2}\,
\frac{z+11z^2+11z^3+z^4}{(z-1)^5}\biggr)\bigg|_{z=e^\lambda},
\end{align*}
where
$$
\lambda=\lambda(C)
=\log2\biggl(1-\frac3{2m}-\frac C{m^2}\biggr)
$$
agrees with our $k/m$ up to $O(m^{-2})$.
Direct computation then shows that
$$
f_m(0)>0.005m^{-2}-100m^{-3}
\quad\text{and}\quad
f_m(0.004)<-0.00015m^{-2}+100m^{-3}
$$
for $m\ge100$. Therefore,
$f_m(0)>110000/m^3$ for $m>2202\cdot10^4$ and
$f_m(0.004)<-110000/m^3$ for $m>734\cdot10^6$,
so that $|f_m(C)|<110000/m^3$ is possible only
if $0<C<0.004$. Comparing this result with~\eqref{E14} we conclude
that, for $k$ and $m>10^9$ satisfying \eqref{eq:1}, we necessarily
have
$$
\frac km=\log2\biggl(1-\frac3{2m}-\frac{C_m}{m^2}\biggr)
$$
with $0<C_m<0.004$.
\end{proof}

%\begin{remark}
%\label{rem:1}
Clearly, the strategy to deduce further terms in the expansion~\eqref{eq:2} remains
the same, but in order to achieve precision $O(m^{-n})$ for an integer $n\ge2$
we have to use the Taylor expansion of $(1-y)^ke^{ky}$ up to $y^{2n+1}$ (each new term
in~\eqref{eq:2} requires two extra terms in the expansion of $(1-y)^ke^{ky}$).
In this way we get
\begin{align}
k
&=cm-\frac{3}{2}c-\biggl(\frac{25}{12}c-3c^2\biggr)m^{-1}
+\biggl(-\frac{73}{8}c+\frac{61}2c^2-25c^3\biggr)m^{-2}
\nonumber\\ &\qquad
+\biggl(-\frac{41299}{720}c+\frac{657}2c^2-598c^3+\frac{1405}4c^4\biggr)m^{-3}+O(m^{-4})\biggr)
\nonumber\\
&\approx 0.69314718m-1.03972077-0.00269758m^{-1}+0.00323260m^{-2}
\nonumber\\ &\qquad
+0.00217182m^{-3}+O(m^{-4}),
\label{km}
\end{align}
where $c=\log2$. However, we do not possess any clear general strategy
to quantify such expansions. Already proving a sharp dependence
on $k$ for the remainder of the $n$-th truncation of the Taylor expansion of
$(1-y)^ke^{ky}$ (like we do for $n=4$ in Lemma~\ref{lem:As}) seems
to be a difficult task. We discuss related problems in Section~\ref{s5}.
%\end{remark}

\begin{proof}[Proof of Corollary~{\rm\ref{cor:1}}]
Let $(m,k)$ be a non-trivial integer solution of~\eqref{EME}. By Moser's result we know that
$m>10^9$. It follows from Theorem~\ref{th:As} that
\begin{equation}
\label{drie}
0<\log 2-\frac{2k}{2m-3} < \frac{0.0111}{(2m-3)^2}.
\end{equation}
By Legendre's theorem, $|\log2-p/q|<1/(2q^2)$ implies that $p/q$ is a convergent
of~$\log2$, while $\log2>p/q$ insures that the index of the convergent is even.
Thus, $2k/(2m-3)$ is a convergent $p_j/q_j$ of the continued fraction of $\log 2$
with $j$~even.
\end{proof}

\section{The proof of the main theorem}
\label{s3}

In this section we prove Theorem~\ref{main}. The restrictions on the prime factorization
of~$q_j$ in that result are established using an argument in the style of Moser given
in the proof of the following lemma.

\begin{Lem}
\label{twee}
Let $(m,k)$ be a solution of {\rm\eqref{EME}} with $k\ge 2$. Let $p$ be a prime
divisor of $2m-3$. If $p-1\mid k$, then
$$
\nu_p(2m-3)=\nu_p(3^{p-1}-1)+\nu_p(k)+1\ge 2.
$$
If $3$~is a primitive root modulo $p$, then $p-1\mid k$.
\end{Lem}

\begin{proof}
Using that $k$ must be even, we find that
\begin{align*}
\sum_{j=1}^{2m-4}j^k
& \equiv \sum_{j=1}^{m-1}j^k+\sum_{j=1}^{m-3}(2m-3-j)^k
\equiv \sum_{j=1}^{m-1}j^k+\sum_{j=1}^{m-3}j^k\pmod{2m-3}
\\
& \equiv m^k+m^k-(m-1)^k-(m-2)^k
\equiv 2(3^k-1)(m-1)^k\pmod{2m-3},
\end{align*}
where we used that $m^k\equiv (2m-3+m)^k\equiv 3^k(m-1)^k\pmod{2m-3}$ and
$(m-2)^k\equiv (2m-3-m+1)^k\equiv (m-1)^k\pmod{2m-3}$.
On applying \eqref{staudt} with $l=2m-3$ and $r=k$ we then obtain that
\begin{equation}
\label{sta}
\frac{2(3^k-1)(m-1)^k}{2m-3}\equiv -\sum_{\substack{p\mid 2m-3\\p-1\mid k}}\frac1p\pmod1.
\end{equation}
If $p\mid 2m-3$ and $p-1\mid k$, the $p$-order of the right-hand side is~$-1$. The $p$-order of
the left-hand side must also be $-1$, that is, we must have
$$
\nu_p(2m-3)=\nu_p(3^k-1)+k\nu_p(m-1)+1=\nu_p(3^{p-1}-1)+\nu_p(k)+1,
$$
where we used that $m-1$ and $2m-3$ are coprime.
Now suppose that $p\mid 2m-3$ and $3$ is a primitive root modulo $p$ (thus $p\mid 3^k-1$ implies
$p-1\mid k$). If $p-1\nmid k$, the $p$-order of the left-hand side is $\le -1$ and $>-1$ on the
right-hand side. Thus, we infer that $p-1\mid k$.
\end{proof}

This completes the required ingredients needed in order to prove the main result.

\begin{proof}[Proof of Theorem {\rm\ref{main}}]
Since by assumption $N\mid k$, we can write $k=Nk_1$ and thus rewrite \eqref{drie} as
\begin{equation}
\label{E}
0<\frac{\log 2}{2N}-\frac{k_1}{2m-3}<\frac{0.0111}{2N(2m-3)^2}.
\end{equation}
We infer that $k_1/(2m-3)=p_j/q_j$ is a convergent to
$(\log 2)/(2N)$ with $j$~even.
Since $p\mid m$ implies $p-1\nmid k$ (see, e.g., Moree \cite[Proposition 9]{Oz}),
we have $(6,q_j)=1$.
We rewrite~\eqref{E} as
$$
0<\frac{\log 2}{2N}-\frac{p_j}{q_j}<\frac{0.0111}{2Nd^2q_j^2},
$$
with $d$ the greatest common divisor of $k_1$ and $2m-3$.
On the other hand,
$$
\frac{\log 2}{2N}-\frac{p_j}{q_j}>\frac1{(a_{j+1}+2)q_j^2},
$$
hence $(a_{j+1}+2)^{-1}<0.0111/(2Nd^2)$, from which the
result follows on also noting that $2m-3\ge q_j$ and invoking Lemma~\ref{twee}
(note that if $\nu_p(q_j)\ge 1$, then $\nu_p(q_j)=\nu_p(2m-3)-\nu_p(k_1)$).
\end{proof}

To prove that $p\mid m$ implies $p-1\nmid k$ one uses that $k$ must be even
and takes $l=m$ in~\eqref{staudt}, showing that $\sum_{p\mid m,\,p-1\mid k}\frac1p$ must
be an integer. Since a sum of reciprocals of distinct primes can never be an integer,
the result follows.

\section{Computation of the continued fractions}
\label{s4}

We make use of conditions \ref{ca}, \ref{cb}, \ref{cc} of Theorem~\ref{main}. We recall that
we expect ${\mathsf E}(\log q_{j(N)}(\xi))\sim c_1N$ for a generic $\xi\in [0,1]$ satisfying
these conditions.
Indeed, on the basis of theoretical results, heuristics and numerical experiments, we conjecture
that $c_1=60\pi^2$.

\begin{table}[ht]%[p]
\begin{center}
\begin{tabular}{|c|r|r|l|c|c|}
\hline
\multicolumn{1}{|c|}{\T $N$} & \multicolumn{1}{|c|}{$j=j(N)$} & \multicolumn{1}{|c|}{$a_{j+1}$} &
\multicolumn{1}{|c|}{$q_j$ \ (rounded down)} & \multicolumn{1}{|c|}{\small $q_j \bmod 6$} & \multicolumn{1}{|c|}{$p = p(q_j)$} \\
\hline
\T $1$ & $642$ & $764$ & $2.383153 \cdot 10^{\,330}$ & $-1$ & $149$ \\
\hline
\T $2$ & $664$ & $1\,529$ & $2.383153 \cdot 10^{\,330}$ & $-1$ & $149$ \\
\hline
\T $2^2$ & $1\,254$ & $21\,966$ & $1.132014 \cdot 10^{\,638}$ & $+1$ & $5$ \\
\hline
\T $2^3$ & $1\,264$ & $43\,933$ & $1.132014 \cdot 10^{\,638}$ & $+1$ & $5$ \\
\hline
\T $2^4$ & $1\,280$ & $87\,866$ & $1.132014 \cdot 10^{\,638}$ & $+1$ & $5$ \\
\hline
\T $2^5$ & $1\,294$ & $175\,733$ & $1.132014 \cdot 10^{\,638}$ & $+1$ & $5$ \\
\hline
\T $2^6$ & $8\,950$ & $26\,416$ & $3.458446 \cdot 10^{\,4\,589}$ & $-1$ & $$ \\
\hline
\T $2^7$ & $8\,926$ & $52\,834$ & $3.458446 \cdot 10^{\,4\,589}$ & $-1$ & $$ \\
\hline
\T $2^8$ & $119\,476$ & $122\,799$ & $1.374540 \cdot 10^{\,61\,317}$ & $+1$ & $$ \\
\hline
\T $2^8 \cdot 3$ & $119\,008$ & $368\,398$ & $1.374540 \cdot 10^{\,61\,317}$ & $+1$ & $$ \\
\hline
\T $2^8 \cdot 3^2$ & $139\,532$ & $782\,152$ & $9.351282 \cdot 10^{\,71\,882}$ & $+1$ & $56\,131$ \\
\hline
\T $2^8 \cdot 3^3$ & $6\,168\,634$ & $1\,540\,283$ & $8.220719 \cdot 10^{\,3\,177\,670}$ & $+1$ & $$ \\
\hline
\T $2^8 \cdot 3^4$ & $22\,383\,618$ & $5\,167\,079$ & $5.128265 \cdot 10^{\,11\,538\,265}$ & $+1$ & $17$ \\
\hline
\T $2^8 \cdot 3^5$ & $155\,830\,946$ & $31\,664\,035$ & $2.257099 \cdot 10^{\,80\,303\,211}$ & $-1$ & $$ \\
\hline
\T $2^8 \cdot 3^5 \cdot 5$ & $351\,661\,538$ & $85\,898\,211$ & $9.729739 \cdot 10^{\,181\,214\,202}$ & $-1$ & $$ \\
\hline
\T $2^8 \cdot 3^5 \cdot 5^2$ & $1\,738\,154\,976$ & $1\,433\,700\,727$ & $1.594940 \cdot 10^{\,895\,721\,905}$ & $+1$ & $5$ \\
\hline
\T $$ & $1\,977\,626\,256$ & $853\,324\,651$ & $1.196828 \cdot 10^{\,1\,019\,133\,881}$ & $-1$ & $$ \\
\hline
\T $2^8 \cdot 3^5 \cdot 5^3$ & $2\,015\,279\,170$ & $4\,388\,327\,617$ & $5.565196 \cdot 10^{\,1\,038\,523\,018}$ & $-1$ & $19$ \\
\hline
\T $$ & $3\,236\,170\,820$ & $2\,307\,115\,390$ & $5.427815 \cdot 10^{\,1\,667\,658\,416}$ & $+1$ & \\
\hline
\T $2^8 \cdot 3^5 \cdot 5^4$ & $2\,015\,385\,392$ & $21\,941\,638\,090$ & $5.565196 \cdot 10^{\,1\,038\,523\,018}$ & $-1$ & $19$ \\
\hline
\T $$ & $3\,236\,257\,942$ & $11\,535\,576\,954$ & $5.427815 \cdot 10^{\,1\,667\,658\,416}$ & $+1$ & $$ \\
\hline
\end{tabular}
\end{center}
\caption{Smallest integers $j$ satisfying conditions \ref{ca}, \ref{cb} and \ref{cc} of Theorem~\ref{main}}
\label{tab:j}
\end{table}

The computation of $(\log 2)/(2N)$ is done in two steps. First, we generate $d$~digits
of~$\log2$. For this we use the $\gamma$-cruncher~\cite{AJY}. With
this program, A.~Yee and R.~Chan computed 31~billion decimal digits
of~$\log2$ in about 24~hours. Second, we set a rational approximation
of $(\log 2)/(2N)$ with a relative error bounded by~$10^{-d}$. Then
partial quotients of the continued fraction of $(\log 2)/(2N)$ are
computed: about $0.97 d$ of them can be evaluated, with safe error
control~\cite{BPR96} (cf.\ the result of Lochs mentioned in Section~\ref{s1}).
We maintain a floating point approximation
of numbers~$q_j$ (rounded down) and residues of $q_j \pmod 6$
by the formula $q_{i+1} = a_{i+1} q_i + q_{i-1}$ for $i \geq 0$,
where $q_0 = 1$ and $q_{-1} = 0$.

Table~\ref{tab:j} was created with the `basic method' of~\cite{BPR96}
for $N \leq 2^8 \cdot 3^4$. It was fast enough to reach the benchmark
$m > 10^{10^7}$ in four days with $50\cdot 10^6$ digits of~$\log2$.
Bit-complexity of this algorithm (or of the indirect or direct
methods~\cite{BPR96}) is quadratic and reaching the $m > 10^{10^{10}}$ milestone
would take centuries.

Some subquadratic GCD algorithms were discovered that have
asymptotic running time $O(n(\log n)^2\log\log n)$~\cite{MOL08}.
A faster version of the program was written: this time a recursive
HGCD method is applied. It is adapted for computing a continued
fraction by using Lemma~3 of~\cite{BPR96} (which is similar to
Algorithm~1.3.13 of~\cite{Cohen93}) for error control. With it
the program leaps over $10^{10^8}$ in just about one hour.
Finally, the new benchmark $m > 10^{10^9}$ is established in no more
than 10~hours with $3\cdot 10^9$ digits of~$\log2$, $N = 1555200$
and condition~\ref{cd}: the first found solution fits
conditions \ref{ca}--\ref{cc}, but not \ref{cd}. With $N = 7776000$, $m > 10^{10^9}$
is achieved for the smallest $j$. See Table~\ref{tab:j}: in the
last column, $p$~is a prime such that $p \in \mathcal{P}(N)$ and
$\nu_p(q_j) = 1$, that is, such that condition~\ref{cd} of Theorem~\ref{main}
is violated.

Now, computation time is not a problem to achieve the $m > 10^{10^{10}}$
milestone, a few days will be sufficient on a computer with a large
amount of memory. We remark that the complexity and
hardware requirement for computation of the digits of~$\log2$, respectively
for computation of its continued fraction expansion, are similar.

\section{Miscellaneous}
\label{s5}

\subsection{{\tt The number of distinct prime factors of $m-1$}}
\label{s5.1}

There is a different application of Theorem~\ref{YVES} suggested
by the work of Brenton and Vasiliu~\cite{BV}, to factorization
properties of the number $m-1$ coming from a non-trivial solution $(m,k)$ of~\eqref{EME}.
A result of Moser~\cite{Moser} (which can also be deduced from the key identity~\eqref{staudt},
cf.\ the proof of Lemma~\ref{twee} above) asserts that
\begin{equation}
\sum_{p\mid m-1}\frac1p+\frac1{m-1}\in\mathbb Z;
\label{pm}
\end{equation}
in particular, the number $m-1$ is square-free. Since the sum of reciprocals
of the first 58~primes is less than~2, we conclude that either
$\omega(m-1)\ge 58$ or the integer
in~\eqref{pm} is equal to~1. In the latter case,
we can apply Curtiss' bound~\cite{Cur} for positive integer solutions
of Kellogg's equation
$$
\sum_{i=1}^n\frac1{x_i}=1,
$$
namely, $\max_i\{x_i\}\le A_n-1$, where the Sylvester sequence
$\{A_n\}_{n\ge1}=\{2,3,7,43,\dots\}$ is defined by the recurrence
$A_n=1+\prod_{i=1}^{n-1}A_i$ (for some further info, see e.g. Odoni \cite{Odoni}).
{}From this result and the estimate
$A_n<(1.066\cdot10^{13})^{2^{n-7}}$, we infer
$$
m<(1.066\cdot10^{13})^{2^{\omega(m-1)-6}},
$$
which together with the lower bound on~$m$ from Theorem~\ref{YVES} yields
$\omega(m-1)\ge 33$.
A similar estimate on the basis of another \eqref{pm}-like identity of Moser implies
that $\omega(m+1)\ge 32$.

\subsection{{\tt Generalized EM equation}}
\label{s5.2}

The method we use in Section~\ref{sec:As} for deriving the asymptotics of $k$ in terms of~$m$
works for the more general equation
\begin{equation}
\label{EMEt}
1^k+2^k+\dots+(m-1)^k=tm^k,
\end{equation}
with $t\in\mathbb N$ fixed, as well. Indeed, the coefficients in the Taylor series expansion
\begin{equation}
\label{TS}
(1-y)^ke^{ky}
=1-\frac k2y^2-\frac k3y^3+\frac{k(k-2)}8y^4+\dotsb
=\sum_{n=0}^\infty g_n(k)y^n
\end{equation}
are polynomials satisfying
\begin{equation}
\label{TSa}
g_0(k)=1, \quad g_1(k)=0,
\qquad\text{and}\qquad
\deg_kg_n(k)=\biggl[\frac n2\biggr], \quad
g_n(0)=0 \quad\text{for $n\ge2$};
\end{equation}
the latter follows from raising the series $(1-y)e^y=1-y^2/2-y^3/3-\dotsb$
to the power~$k$. In these settings, equation~\eqref{EMEt} becomes
\begin{align*}
t
&=\sum_{j=1}^{m-1}\biggl(1-\frac jm\biggr)^k
=\sum_{j=1}^{m-1}e^{-kj/m}\sum_{n=0}^\infty g_n(k)\biggl(\frac jm\biggr)^n
\\
&=\sum_{n=0}^\infty\frac{g_n(k)}{m^n}\sum_{j=1}^{m-1}j^ne^{-jk/m}
\\ \intertext{(since $\sum_{j=m}^\infty j^ne^{-jk/m}=O(m^ne^{-k})$)}
&\sim\sum_{n=0}^\infty\frac{g_n(k)}{m^n}\sum_{j=1}^\infty j^ne^{-jk/m}
=\sum_{n=0}^\infty\frac{g_n(k)}{m^n}
\biggl(\biggl(z\frac{\d}{\d z}\biggr)^n\frac z{1-z}\biggr)\bigg|_{z=e^{-k/m}}
\\
&=\sum_{n=0}^\infty\frac{g_n(k)}{m^n}(-1)^n
\biggl(\biggl(z\frac{\d}{\d z}\biggr)^n\frac1{z-1}\biggr)\bigg|_{z=e^{k/m}},
\end{align*}
hence in the notation $\lambda=k/m$ and $x=1/m$ we have
\begin{equation}
\label{ast}
t=\sum_{n=0}^\infty g_n\biggl(\frac\lambda x\biggr)(-x)^n
\biggl(\biggl(z\frac{\d}{\d z}\biggr)^n\frac1{z-1}\biggr)\bigg|_{z=e^\lambda}.
\end{equation}
Searching $\lambda$ in the form $\lambda=c_0+c_1x+c_2x^2+\dotsb$,
we find successively
\begin{gather*}
c_0=c(t)=\log\biggl(1+\frac1t\biggr)=\log\frac{t+1}t,
\qquad c_1=-\biggl(t+\frac12\biggr)c,
\\
c_2=\biggl(t+\frac12\biggr)^3c^2-\biggl(t+\frac12\biggr)^2c
-\frac14\biggl(t+\frac12\biggr)c^2+\frac c6,
\end{gather*}
and so on. Note that $c_n(-(t+1))=(-1)^{n+1}c_n(t)$ for $n=0,1,2,\dots$;
this reflects the equivalence of equation \eqref{EMEt} and
\begin{equation}
\label{EMEtdual}
1^k+2^k+\dots+(m-1)^k+m^k=(t+1)m^k.
\end{equation}

{}From this asymptotics we see that
\begin{align}
\frac{2k}{2m-t_1}
&=c+\frac{{t_1^3}c^2-2{t_1^2}c-t_1c^2+4c/3}{2(2m-t_1)^2}+O\biggl(\frac1{(2m-t_1)^3}\biggr),
\label{CFt}
\end{align}
where $t_1=2t+1$ and $c=\log(1+1/t)$.
It can be checked that for all positive integers $t$ we have the inequality
$$
-0.22<{t_1^3}c^2-2{t_1^2}c-t_1c^2+\frac{4c}3<0,
$$
and hence $2k/(2m-2t-1)$ is a convergent (with even index) of this logarithm
$c=\log(1+1/t)$ for $m$~large enough.

\subsection{{\tt Saddle-point method}}
\label{s5.3}

A different approach to treat the asymptotic behaviour of $k$ in terms of $m$
for $k$ and $m$ satisfying \eqref{EME} (or, more generally, \eqref{EMEt})
is based on the integral representation
$$
1^k+2^k+\dots+(m-1)^k
=\frac{\Gamma(k)}{2\pi i}\int_{C-i\infty}^{C+i\infty}
\frac{e^{mz}}{(e^z-1)z^{k+1}}\,\d z,
$$
where $C$~is an arbitrary positive real number (cf.~\cite[p.~273]{Delange}).
On noting that
$$
\frac{e^{mz}}{e^z-1}
=\frac{e^{(m-1)z}}{1-e^{-C}}\biggl(1+\frac{1-e^{z-C}}{e^z-1}\biggr)
$$
one obtains, on taking $C=(k+1)/(m-1)$ and after invoking some rather
trivial estimates, that
\begin{equation}
\label{hubert}
1^k+2^k+\dots+(m-1)^k=\frac{(m-1)^k}{1-e^{-(k+1)/(m-1)}}\bigl(1+\rho_k(m)\bigr),
\end{equation}
with
$$
|\rho_k(m)|<\frac{\sqrt{2(k+1)}C}{\sqrt\pi(k-1)(e^C-1)}.
$$
(This part of the argument is due to Delange; for more details see
\cite[pp.~273--274]{Delange}.) By~\eqref{ineq},
$C$~is bounded and we infer that $|\rho_k(m)|=O(k^{-1/2})=O(m^{-1/2})$.
On putting $m^k$ on the left-hand side of~\eqref{hubert}
and using $(1-1/m)^m=\exp(-1+O(m^{-1}))$, we immediately conclude
that, as $m\to\infty$,
$$
\frac km=\log2+O\biggl(\frac1{\sqrt{m}}\biggr),
$$
where the implied constant is absolute. A more elaborate analysis,
using the saddle-point method, will very likely allow one as many terms
in the latter expansion as required.

\subsection{{\tt Experimental asymptotics}}
\label{s5.4}

It is worth mentioning a fast experimental approach of doing asymptotics
like~\eqref{km}. Given numerically a few hundred terms of a sequence
$s=\{s_n\}_{n\ge1}$ that one believes has an asymptotic expansion
in inverse powers of $n$, one can try to apply the $\mathtt{asymp}_k$
trick, a simple but often powerful method to numerically determine the
coefficients in the ansatz
$$
s_n\sim c_0+\frac{c_1}{n}+\frac{c_2}{n^2}+\dotsb.
$$
As a second step one tries to identify the so-found coefficients
with (linear combinations of) known constants. Thus, one arrives
at a conjecture that hopefully can be turned into a proof.
For more details and some `victories' achieved by
the $\mathtt{asymp}_k$ method, see Gr\"unberg and Moree \cite{GM}.

D.~Zagier has applied this trick to the sequence of $k=k(m)$ obtained from
\eqref{EME} on letting $m$ run through the first thousand values.
Excellent agreement with our theoretical results was obtained in this way.

\medskip
{\tt Acknowledgements}. The second author is very indebted to Jerzy Urbanowicz
for involving him in the early nineties in his EM research (with \cite{MRU} as
visible outcome). The second and third author would like to thank D.~Zagier
for verifying some of our results using the $\mathtt{asymp}_k$ trick and for some
informative discussions regarding the saddle-point method (reflected in the final section). H.~te Riele provided us
with the unpublished report
\cite{Best}, which became the `initial spark' for the current project. C.~Baxa pointed out the
relevance of \cite{Harman} to us. Further thanks are due to T.~Agoh and I.~Shparlinski.

This research was carried out whilst the third author was visiting in the Max Planck Institute for Mathematics
(MPIM) and the Hausdorff Center for Mathematics (HCM) financially supported by these institutions.
He and the second author thank the MPIM and HCM for providing such a nice research environment.


\begin{thebibliography}{99}

\bibitem{Alzer}
\textsc{H.~Alzer},
\"Uber eine Verallgemeinerung der Bernoullischen Ungleichung,
\emph{Elem. Math.} \textbf{45} (1990), 53--54.

\bibitem{BB}
\textsc{J.~Borwein} and \textsc{D.~Bailey},
\emph{Mathematics by experiment. Plausible reasoning in the 21st century}, 2nd edition
(A K Peters, Ltd., Wellesley, MA, 2008).

\bibitem{Best}
\textsc{M.\,R.~Best} and \textsc{H.\,J.\,J.~te Riele},
On a conjecture of Erd\H{o}s concerning sums of powers of integers,
Report NW 23/76 (Mathematisch Centrum Amsterdam, 1976).

\bibitem{BPR96}
\textsc{R.\,P.~Brent}, \textsc{A.\,J.~van der Poorten}, and \textsc{H.~te Riele},
A comparative study of algorithms for computing continued fractions of algebraic numbers,
\emph{Algorithmic number theory} (Talence, 1996),
Lecture Notes in Comput. Sci. \textbf{1122} (Springer, Berlin, 1996), 35--47.

\bibitem{BV}
\textsc{L.~Brenton} and \textsc{A.~Vasiliu},
Znam's problem,
\emph{Math. Mag.} \textbf{75} (2002), 3--11.

\bibitem{graphs}
\textsc{W.~Butske}, \textsc{L.\,M.~Jaje} and \textsc{D.\,R.~Mayernik},
On the equation $\sum_{p\mid N}\frac1p+\frac1N=1$, pseudoperfect numbers, and perfectly weighted graphs,
\emph{Math. Comp.} \textbf{69} (2000), 407--420.

\bibitem{Carlitz}
\textsc{L.~Carlitz},
The Staudt--Clausen theorem,
\emph{Math. Mag.} \textbf{34} (1960/1961), 131--146.

\bibitem{Cohen93}
\textsc{H.~Cohen},
\emph{A course in computational algebraic number theory},
Graduate Texts in Math. \textbf{138} (Springer, Berlin, 1993).

\bibitem{Cur}
\textsc{D.\,R.~Curtiss},
On Kellogg's Diophantine problem,
\emph{Amer. Math. Monthly} \textbf{29} (1922), 380--387.

\bibitem{Delange}
\textsc{H.~Delange},
Sur les z\'eros r\'eels des polyn\^omes de Bernoulli,
\emph{Ann. Inst. Fourier \textup(Grenoble\textup)} \textbf{41} (1991), 267--309.

\bibitem{E}
\textsc{P.~Erd\H{o}s},
Advanced Problem 4347,
\emph{Amer. Math. Monthly} \textbf{56} (1949), 343.

\bibitem{Gerber}
\textsc{L.~Gerber},
An extension of Bernoulli's inequality,
\emph{Amer. Math. Monthly} \textbf{75} (1968), 875--876.

%\bibitem{GS} X.~Gourdon and P.~Sebah, The logarithm constant: $\log 2$,\\
%\texttt{http://numbers.computation.free.fr/Constants/Log2/log2.html}

%\bibitem{Galambos}
%\textsc{J.~Galambos},
%The distribution of the largest coefficient in continued fraction expansions,
%\emph{Quart. J. Math. Oxford Ser.} (2)  \textbf{23}  (1972), 147--151.

\bibitem{GM}
\textsc{D.\,B.~Gr\"unberg} and \textsc{P.~Moree},
Sequences of enumerative geometry: congruences and asymptotics
(with an appendix by D.~Zagier),
\emph{Experiment. Math.}  \textbf{17}  (2008), 409--426.

\bibitem{Guy}
\textsc{R.\,K.~Guy},
\emph{Unsolved problems in number theory}, 3rd edition,
Problem Books in Mathematics (Springer, New York, 2004).

\bibitem{Harman}
\textsc{G.~Harman} and \textsc{K.\,C.~Wong},
A note on the metrical theory of continued fractions,
\emph{Amer. Math. Monthly } \textbf{107}  (2000),  834--837.

\bibitem{Hooley}
\textsc{C.~Hooley},
On Artin's conjecture,
\emph{J. Reine Angew. Math.} \textbf{225} (1967), 209--220.

\bibitem{JL}
\textsc{H.~Jager} and \textsc{P.~Liardet},
Distributions arithm\'etiques des d\`enominateurs de convergents de fractions continues,
\emph{Nederl. Akad. Wetensch. Indag. Math.} \textbf{50} (1988), 181--197.

\bibitem{Kanada}
\textsc{Y.~Kanada},
Kanada $\pi$-Laboratory,
available at \texttt{http://\allowbreak www.\allowbreak super-computing.\allowbreak org/}.

\bibitem{Kellner}
\textsc{B.\,C.~Kellner},
\"Uber irregul\"are Paare h\"ohere Ordnungen,
Diplomarbeit
(Mathematisches Institut der Georg-August-Universit\"at zu G\"ottingen,
Germany, 2002);
available at \texttt{http://\allowbreak www.bernoulli.org/\allowbreak\~{}bk/irrpairord.pdf}.

\bibitem{K}
\textsc{B.~Krzysztofek},
The equation $1^n+\ldots+m^n=(m+1)^nk$,
\emph{Wyz. Szkol. Ped. w. Katowicech-Zeszyty Nauk. Sekc. Math.} \textbf{5} (1966), 47--54. (Polish)

\bibitem{L}
\textsc{P.~L\'evy},
Sur le d\'eveloppement en fraction continue,
\emph{Compositio Math.} \textbf{3} (1936), 286--303.

\bibitem{Liardet}
\textsc{P.~Liardet}
Propri\'et\'es arithm\'etiques presque s\^ures des convergents,
\emph{S\'eminaire de th\'eorie des Nombres de Bordeaux} (1986-87), exp. no.~36, 20~pp.

%\bibitem{MB}
%\textsc{S.\,J.~Miller} and \textsc{R.~Takloo-Bighash},
%\emph{An invitation to modern number theory}
%(Princeton University Press, Princeton, NJ, 2006).

\bibitem{Moeckel}
\textsc{R.~Moeckel},
Geodesics on modular surfaces and continued fractions,
\emph{Ergodic Theory Dynamical Systems} \textbf{2} (1982), 69--83.

\bibitem{MOL08}
\textsc{N.~M\"oller},
On Sch\"onhage's algorithm and subquadratic integer GCD computation,
\emph{Math. Comp.} \textbf{77} (2008), 589--607.

\bibitem{Canada}
\textsc{P.~Moree},
On a theorem of Carlitz--von Staudt,
\emph{C.~R. Math. Rep. Acad. Sci. Canada} \textbf{16} (1994),  166--170.

\bibitem{Oz}
\textsc{P.~Moree},
Diophantine equations of Erd\H{o}s--Moser type,
\emph{Bull. Austral. Math. Soc.} \textbf{53} (1996), 281--292.

%\bibitem{MC}
%\textsc{P.~Moree} and \textsc{J.~Cazaran},
%On a claim of Ramanujan in his first letter to Hardy.
%\emph{Exposition. Math.}  \textbf{17}  (1999),  289--311.

\bibitem{MRU}
\textsc{P.~Moree}, \textsc{H.~te Riele}, and \textsc{J.~Urbanowicz},
Divisibility properties of integers $x$, $k$ satisfying $1^k+\dots+(x-1)^k=x^k$,
\emph{Math. Comp.} \textbf{63} (1994), 799--815.

\bibitem{Moser}
\textsc{L.~Moser},
On the diophantine equation $1^n+2^n+3^n+\dots+(m-1)^n=m^n$,
\emph{Scripta Math.} \textbf{19} (1953), 84--88.

\bibitem{Odoni}
\textsc{R.\,W.\,K.~Odoni},
On the prime divisors of the sequence $w_{n+1}=1+w_1\cdots w_n$,
\emph{J. London Math. Soc.} (2) \textbf{32} (1985), 1--11.

\bibitem{Urbi}
\textsc{J.~Urbanowicz},
Remarks on the equation $1^k+2^k+\dots+(x-1)^k=x^k$,
\emph{Indag. Math.} \textbf{50} (1988), 343--348.

\bibitem{AJY}
\textsc{A.\,J.~Yee},
$\gamma$-cruncher\,---\,A Multi-Threaded Pi-Program,
available at \texttt{http://\allowbreak www.\allowbreak numberworld.\allowbreak org/}.

\end{thebibliography}
\end{document}